\documentclass[11pt,reqno]{amsart}
\usepackage{graphicx}
\usepackage{amsmath,amssymb}
\usepackage{latexsym}
\usepackage{mathrsfs}
\usepackage{color}
\usepackage[utf8]{inputenc}
\usepackage{amsmath,amssymb,amsfonts,amsthm,mathrsfs,enumerate,cite,xcolor,graphicx}
\usepackage[pdfpagelabels,linktocpage]{hyperref}
\usepackage{epstopdf}
\usepackage{graphicx}
\usepackage{amssymb,amsmath}
\usepackage{latexsym}
\usepackage{epsfig}

    \newcommand{\Z}{{\mathbb{Z}}}
\newcommand{\Y}{{\mathbb{Y}}}

    \newcommand{\R}{\mathbb{R}}
    \newcommand{\be}{\begin{eqnarray}}
    \newcommand{\ee}{\end{eqnarray}}

\newcommand{\cof}{{\rm cof}}
\renewcommand{\det}{{\rm det}}
\newcommand{\D}{\mathcal{D}}

\renewcommand{\O}{\Omega}

\newcommand{\md}{\mathrm{d}}

\newcommand{\I}{{\mathcal{I}}}

\usepackage{xargs}

\numberwithin{equation}{section}
\mathchardef\emptyset="001F

\newtheorem{theorem}{Theorem}[section]
\newtheorem{proposition}[theorem]{Proposition}

\newtheorem{remark}[theorem]{Remark}
\newtheorem{definition}[theorem]{Definition}

\begin{document}
\title{Elastoplastic deformations of  layered structures}

\author[Daria Drozdenko]{Daria Drozdenko}
\address[D. Drozdenko]{Faculty of Mathematics and Physics, Charles University, Ke Karlovu 5, CZ-121 16, Prague 2, Czechia}
\email{Daria.Drozdenko@mff.cuni.cz}

\author[Michal Knapek]{Michal Knapek}
\address[M. Knapek]{Faculty of Mathematics and Physics, Charles University, Ke Karlovu 5, CZ-121 16, Prague 2, Czechia}
\email{knapek@karlov.mff.cuni.cz}

\author[M. Kru\v{z}\'{i}k]
{Martin Kru\v{z}\'{i}k}
\address[M. Kru\v{z}\'{i}k]{Institute of Information Theory and
  Automation, Czech Academy of Sciences,
Pod vod\'arenskou ve\v z\'\i\ 4, CZ-182 00, Prague 8, Czechia and 
Faculty of Civil Engineering, Czech Technical University, Th\'{a}kurova 7, CZ-166 29, Prague 6, Czechia}
\email{kruzik@utia.cas.cz}

\author[Kristi\'{a}n  M\'{a}this]{Kristi\'{a}n M\'{a}this}
\address[K. M\'{a}this]{Faculty of Mathematics and Physics, Charles University, Ke Karlovu 5, CZ-121 16, Prague 2, Czechia}
\email{Kristian.Mathis@mff.cuni.cz}

\author[Karel \v{S}vadlenka]{Karel \v{S}vadlenka}
\address[K. \v{S}vadlenka]{Kyoto University, Graduate School of Science, Kitashirakawa Oiwake-cho, Sakyo-ku, Kyoto, 606-8502 Japan}
\email{karel@math.kyoto-u.ac.jp} 

\author[Jan Valdman]{Jan Valdman}
\address[Jan Valdman]{Institute of Information Theory and
  Automation, Czech Academy of Sciences,
Pod vod\'arenskou ve\v z\'\i\ 4, CZ-182 00, Prague 8, Czechia and Faculty of Science,
University of South Bohemia, Brani\v{s}ovsk\'{a} 1645/31a,
CZ-370 05 \v{C}esk\'{e}~Bud\v{e}jovice, Czechia}
\email{jan.valdman@utia.cas.cz} 

\maketitle

\section*{Abstract}
We formulate a large-strain model of single-slip crystal elastoplasticity in the framework of energetic solutions.  Numerical performance of the model is compared with lab experiments on compression of a stack of note papers. 
\section{Introduction}
Elastoplasticity at large strains is an area of ongoing research that brings together contributions from modeling, mathematical, analysis, numerical simulations, and mechanical experiments. For the mathematical analysis of elastoplastic models, it is often convenient to use powerful tools from the calculus of variations, which are now able to treat quasistatic evolutionary rate-independent problems as well, see, e.g. \cite{cckham, ortiz-repetto}, or  \cite{mielke1}. The existence of solutions could be ensured by assuming generalized convexity of the strain energy, such as polyconvexity \cite{ball77}  but more  general material behavior may contradict this assumption. For example, this is manifested in shape-memory alloys (SMA) \cite{wlsc}, some magnetostrictive \cite{desimone} and ferroelectric materials \cite{shu}.

As a remedy, one can then recourse to higher-gradient regularizations, where the stored energy density $W$ also depends, e.g., on the second gradient of the deformation. From a mathematical point of view, this adds compactness to the model, which is instrumental in proving the existence of solutions by the direct method \cite{ball81}. Materials with such constitutive equations are referred to as non-simple and were introduced by Toupin \cite{toupin, toupin2}. Since then, many authors have elaborated on the concept so that its thermodynamical aspects  are also better understood, cf.~\cite{capriz, hooke2, forest, forest2, greenRivlin,  mierouNum, podio}. We will, however, restrict ourselves to material models where polyconvexity is sufficient for the existence of a solution. Extension to more general material models can be found in many works, here we mention \cite{RIS}.
In what follows, we will deal with a model that, in addition to elasticity, also includes the plastic behavior of the material. Plastic strain acts as an internal variable that influences the elastic behavior of the material. 
In the large strain setting, we assume that 
the deformation gradient $F$ is decomposed by means of the Kr\"{o}ner-Lee decomposition as 
\begin{align}
    F=F^{\rm e}F^{\rm p} \ ,
\end{align}
where $F^{\rm p}:\Omega\to\R^{d\times d}$ and $\det F^{\rm p}=1$. We will set $P=(F^{\rm p})^{-1}$. Here $\Omega \subset\R^d$ is a bounded Lipschitz domain representing the specimen.
Plastic strain, $F^{\rm p}$,  represents disarrangements of atoms and $F^{\rm e}$ describes stretching and a rotation of the lattice.  Unlike $F$, $F^{\rm e}$ or $F^{\rm p}$ do not correspond to deformations, that is, they do not necessarily have to be curl-free, in general. However, we refer to a dislocation-free setting of the problem treated in \cite{mkdmus}.
Here we assume that $F^{\rm p}$ corresponds to the single-slip plastic strain \cite{gurtin}, which is defined as 
$$F^{\rm p}=I+\gamma s\otimes  m \ ,$$ 
where $I\in\R^{d\times d}$ is the identity matrix, $s,m\in\R^d$ are mutually perpendicular unit vectors defining the glide direction $s$ and the slip plane normal $m$. Finally,  $\gamma:\Omega\to\R$ denotes the slip or microshear-strain and it measures the amount of plastic strain. Note that if $\gamma=0$ then $F^{\rm p}=I$, that is, $F=F^{\rm e }$ and the deformation is purely elastic.
It is also easy to see that  $\det (I+\gamma s\otimes m)= 1+\gamma s\cdot m=1$  and $P=I-\gamma s\otimes m$. We refer, e.g., to \cite{grandi1,grandi2, cckham,critplast,jkmk} for recent works on elastoplasticity.


The elastic behavior of the material is described by a stored energy density  $$W:\Omega\times\R^{d\times d}\to\R\cup\{+\infty\},$$  such that the first Piola-Kirchhoff stress tensor $S$ 
is defined for almost every $x\in\O$ as 
\begin{align}
    S(x)=\frac{\partial W(x,F^{e}(x))}{\partial F}\ .
\end{align}
We assume that $W$ is polyconvex \cite{ball77,ciarlet,dacorogna}, that is, there exists a convex and lower semicontinuous Carath\'{e}odory function  $h(x,\cdot):\R^{d\times d}\times\R^{d \times d}\times\R\to\R\cup\{+\infty\}$ such that for almost every $x\in\Omega$
\begin{align}
    W(x,F)= h(x,F,\cof F,\det F) \ \text{ for every } F\in\R^{d \times d} \ .
\end{align}
Moreover, we require that $W$ is frame-indifferent, coercive, and penalizes extreme compression and a change of orientation, i.e., it holds
\begin{subequations}\label{ass-E-1}
\begin{align}
&W(x,F)=W(x,RF) \quad  \text{ for every } R\in\text{SO}(d) \ , \\
&W(x,F)\ge C|F|^p-c_0 \quad \text{ for some } C, c_0>0, p>1 \ , \\
&W(x,F)\to +\infty \quad \text{ if } \det F\to 0 \text{ and } W(x,F)=+\infty \text{ if } \det\, F\le 0\ 
\end{align}
\end{subequations}
and for every $F\in\R^{d\times d}$ and all $x\in\Omega$.
Furthermore, we assume that the energy stored in dislocations is described by a Carath\'{e}odory function $w:\Omega\times \R\to\R$
such that 
\begin{align}\label{ass-gamma}
    w(x,\gamma)\ge C|\gamma|^r-c_0 \text{ for some } C,c_0>0, r>1, \text{ and  all } x\in\Omega \text{ and all } \gamma\in\R.
\end{align}
The functional
\begin{align} \label{energy_gen} \mathcal{I}(t,y,\gamma)&:=\int_\O W(x,\nabla y(x)(I-\gamma(x)s\otimes m))\,\md x+\int_\Omega w(x,\gamma(x))\,\md x\nonumber\\
&+\epsilon\int_\O |s\otimes m\otimes\nabla \gamma(x)|^\alpha\,\md x-L(t,y) \ ,\end{align}
expresses the potential energy in our system. The third term with $\epsilon>0$ depending on $\nabla\gamma$ is the {\it  plastic-strain-gradient} energy penalizing spatial variations of $\gamma$, see, e.g.~\cite{gurtin}. The last term expresses the work done by external force densities $f$ and $g$:  
\be\label{loading1}
L(t,y):=\int_\O f(t)\cdot y\,\md x +\int_{\Gamma_1} g(t)\cdot y\,\md S\ .\
\ee
Here $\Gamma_1\subset\partial\O$ is a part of the boundary where we prescribe some traction.
The appearance of plastic deformation is related to energy dissipation depending on the rate of change of $\gamma$, that is, on $\dot\gamma$. Here, we follow \cite[Formula (54)]{gurtin} where the specific dissipation is given by 
\begin{align}
    \delta(\dot\gamma)=\sigma|\dot\gamma|\ ,
\end{align}
where $\sigma:\Omega\to[\sigma_0,+\infty)$  is the so-called {\it slip resistance} with $\sigma_0>0$.
Consequently, the global dissipation between two states $\gamma_1$ and $\gamma_2$ is defined as 
\begin{align}
    \D(\gamma_1,\gamma_2)=\int_\Omega \sigma(x)|\gamma_1(x)-\gamma_2(x)|\,{\rm d}x\ .
\end{align}

\subsection{Energetic solution}
In order to find a quasistatic evolution of the system, Mielke, Theil, and Levitas \cite{AMFTVL} came up with the following definition of the energetic solution, which conveniently overcomes the non-smoothness of dissipation and is generally very flexible.  Moreover, it fully exploits the possible variational structure of the problem and allows for   a very wide class of 
energy and dissipation functionals. This concept has versatile applications to many problems in the continuum mechanics of solids. Additionally,  working with $\I$ and $\D$ directly enables us to  include higher-order  gradients of $y$ in the model or to require the integrability of some functions of $\nabla y$ if needed.

Let $\Y$ and $\Z$  be the sets of admissible deformations and slips (usually subsets of a Sobolev space) and suppose that the evolution of $y(t)\in \Y$ and $\gamma(t)\in \Z$ is studied during a time interval $[0,T]$ for the time horizon $T>0$.
The following two properties characterize the energetic solution: \\
\noindent
(i) Stability inequality - $\forall t\in[0,T],\, \tilde z\in \Z,\, \tilde{y}\in \Y$:
\be\label{stblt}\mathcal{I}(t,y(t),z(t))\le \mathcal{I}(t,\tilde y,\tilde z)+\mathcal{D}(z(t),\tilde z)\ee

\noindent
(ii) Energy balance - $\forall\ 0\le t\le T$:
\be
\label{e-bal}
\mathcal{I}(t,y(t),\gamma(t))+{\rm Var}(\mathcal{D},\gamma;[0,t])  = \mathcal{I}(0,y(0),\gamma(0)) +\int_0^t \dot L(\xi,y(\xi))\,\md \xi \ ,     
\ee
$$\text{where }
{\rm Var}(\mathcal{D},\gamma;[s,t]):=\sup\left\{\sum_{i=1}^N \mathcal{D}(\gamma(t_i),\gamma(t_{i-1}));\ \{t_i\} \mbox{ partition of } [s,t]\right\}.$$  
  
\begin{definition}\label{en-so}   
The mapping $t\mapsto(y(t),\gamma(t))\in \Y\times \Z $ is an energetic solution to the problem
$(\mathcal{I},\mathcal{D}, L)$ if the stability inequality and the energy balance 
are satisfied for all $t\in [0,T]$.   
\end{definition}

To prove the existence of an energetic solution, we follow the standard strategy described in, e.g. \cite{mami2, RIS}, where this program is specialized in elastoplasticity. 

Let $\Gamma_0\subset\partial\O\setminus\Gamma_1$ be of positive $d-1$ dimensional Lebesgue measure.
Let $p_y>d$ and let  be such that $\int_\O W(x,\nabla y_0(x))\,\md x<+\infty$. Define 
\begin{align*}
&\Y:=\{y\in W^{1,p_y}(\Omega;\R^d): y=y_0 \text{ on } \Gamma_0\} \ , 
\qquad
\Z:=\{\gamma\in L^{r}(\Omega)\cap W^{1,\alpha}(\Omega)\} \ .
\end{align*}
\begin{proposition}\label{prop:existence}
Assume that $\mathcal{I}$ and $\mathcal{D}$ are as above, $L\in C^1([0,T]; W^{1,p_y}(\O;\R^d)^*)$, \eqref{ass-E-1} and \eqref{ass-gamma} hold, $1/p+1/r=1/p_y<1/d$, and $\alpha>1$. Let the initial condition $(y^0,\gamma^0)\in\Y\times\Z$ be stable.  Then an energetic solution exists.
\end{proposition}

{\it Sketch of  proof.} The proof can be obtained following the one in \cite{mami2} or \cite{fm}, where a more general setting is considered. Here, we only sketch it in a few steps for the reader's convenience and for the case of time-independent boundary condition $y_0$.

\noindent
Step 1: Consider a partition $$0=t^0_\tau<t^1_\tau<\ldots<t^K_\tau=T \ ,$$ set  $\tau=\max_i(t_i-t_{i-1})$ and suppose that the partition for $N+1$ is a refinement of the partition with $N$ time steps.  Take the initial condition $(y^0_\tau,\gamma^0_\tau)=(y^0,\gamma^0)\in\Y\times\Z$.
Define the following sequence of minimization problems:
For $k=1,\ldots, K$    solve 
\begin{align}\label{incremental}
\min_{(y,\gamma)\in\Y\times \Z}\mathcal{I}(t_\tau^k,y,\gamma)+\mathcal{D}(\gamma_\tau^{k-1},\gamma)\ 
\end{align}
and denote a solution by $(y_\tau^k,\gamma_\tau^k)$. The existence of a solution is a standard application of the direct method of calculus of variations.

\noindent 
Step 2: The solutions to \eqref{incremental} are stable. Moreover, we have the following:
 \begin{align}\label{e-disc}
\int_{t_\tau^{k-1}}^{t_\tau^k}\partial_t \mathcal{I}(s,y^k_\tau, \gamma^k_\tau)\,\md s \le \mathcal{I}(t^k_\tau,y^k_\tau, \gamma^k_\tau)+\mathcal{D}(\gamma_\tau^{k-1},\gamma_\tau^k)-\mathcal{I}(t^{k-1}_\tau,y^{k-1}_\tau, \gamma^{k-1}_\tau)\\ \nonumber
\le \int_{t_\tau^{k-1}}^{t_\tau^k}\partial_t \mathcal{I}(s,y^{k-1}_\tau,\gamma^{k-1}_\tau)\,\md s\ .
\end{align}
Take $(\tilde y,\tilde{\gamma})\in\Y\times\Z$. We have $\mathcal{I}(t^k_\tau,y^k_\tau,\gamma^k_\tau)+\mathcal{D}(\gamma^{k-1}_\tau,\gamma^k_\tau)\le \mathcal{I}(t^k_\tau,\tilde y,\tilde{\gamma})+\mathcal{D}(\tilde{\gamma}, \gamma^{k-1}_\tau)$. We further estimate $\mathcal{D}(\tilde{\gamma},\gamma_\tau^{k-1})-\mathcal{D}(\gamma_\tau^{k-1},\gamma_\tau^{k})\le \mathcal{D}(\tilde{\gamma}, \gamma^{k}_\tau)$, which proves the stability. 
The upper estimate in (\ref{e-disc}) follows by checking the minimality of $(y_\tau^k,\gamma_\tau^k)$ against $(y_\tau^{k-1},\gamma_\tau^{k-1})$, that is, 
\begin{align*}
\mathcal{I}(t^k_\tau,y^k_\tau,\gamma^k_\tau)+\mathcal{D}(\gamma_\tau^k,\gamma^{k-1}_\tau)&\le \mathcal{I}(t^k_\tau,y^{k-1}_\tau,\gamma^{k-1}_\tau)\\
&=\mathcal{I}(t^{k-1}_\tau,y^{k-1}_\tau,\gamma^{k-1}_\tau )+\int_{t_\tau^{k-1}}^{t_\tau^k}\partial_t \mathcal{I}(s,y^{k-1}_\tau, \gamma^{k-1}_\tau)\,\md s\ .\end{align*}

The lower estimate in (\ref{e-disc}) is implied by the stability of $(y_\tau^{k-1},\gamma_\tau^{k-1})$, that is,
 \begin{align}\mathcal{I}(t^{k-1}_\tau,y^{k-1}_\tau, \gamma^{k-1}_\tau)\le \mathcal{I}(t^{k-1}_\tau,y^{k}_\tau, \gamma^k_\tau)+\mathcal{D}(\gamma_\tau^{k-1},\gamma_\tau^k)\nonumber\\
 =\mathcal{I}(t^{k}_\tau,y^{k}_\tau, \gamma^k_\tau)+\mathcal{D}(\gamma_\tau^{k-1},\gamma_\tau^k)-\int_{t_\tau^{k-1}}^{t_\tau^k}\partial_t \mathcal{I}(s,y^k_\tau, \gamma^k_\tau)\,\md s\ .\end{align}

Then define the piecewise constant interpolants constructed from $\{(y^k,\gamma^k)\}_k$ and denote them by $(y_\tau^K,\gamma_\tau^K)$. In particular, we define 
\be (y^K_\tau(t),\gamma^K_\tau(t)):=(y^{k-1}_\tau, \gamma^{k-1}_\tau) \mbox{ if } t\in[t^{k-1}_\tau, t_\tau^k)\ ,\  (y_\tau(T), \gamma_\tau(T)):=(y^K_\tau,\gamma^K_\tau)\ .\ee 
Using  \eqref{ass-E-1} and \eqref{ass-gamma}  we get the following apriori bounds which are independent of $\tau$:
$$
\|y^K_\tau\|_{L^\infty((0,T);W^{1,p_y}(\O;\R^d))}\le C\ ,
$$
$$
\|\gamma^K_\tau\|_{L^\infty((0,T);L^r(\O))}\le C\ ,
$$
and 
$$
{\rm Var}(\mathcal{D},\gamma^K_\tau;[0,T])\le C \ .$$

\noindent
Step 3: The existence of an energetic solution is now obtained by passing to the limit as $K\to\infty$ (as the time discretization is refined) in the energy inequality proved in Step 2 and checking the stability of the limit. Note that the dissipation functional $\mathcal{D}:\Z\times \Z\to\R$ is sequentially continuous with respect to the weak $W^{1,\alpha}(\O)$ topology.
\hfill
$\Box$

\begin{remark}
It is shown in \cite{mami2} that time-dependent Dirichlet boundary conditions can be considered in the above proposition  if 
$y_0\in  C^1([0,T]\times\O;\R^d)$, $\nabla y_0\in BC^1([0,T]\times\O;\R^{d\times d})$ and $|(\nabla y_0)^{-1}|\in L^\infty(\O)$. Here $BC^1$ stands for bounded and continuously differentiable maps.
\end{remark}

\section{Compression experiment}
\begin{figure}[h]
\begin{center}
\includegraphics[width = 0.7\textwidth]{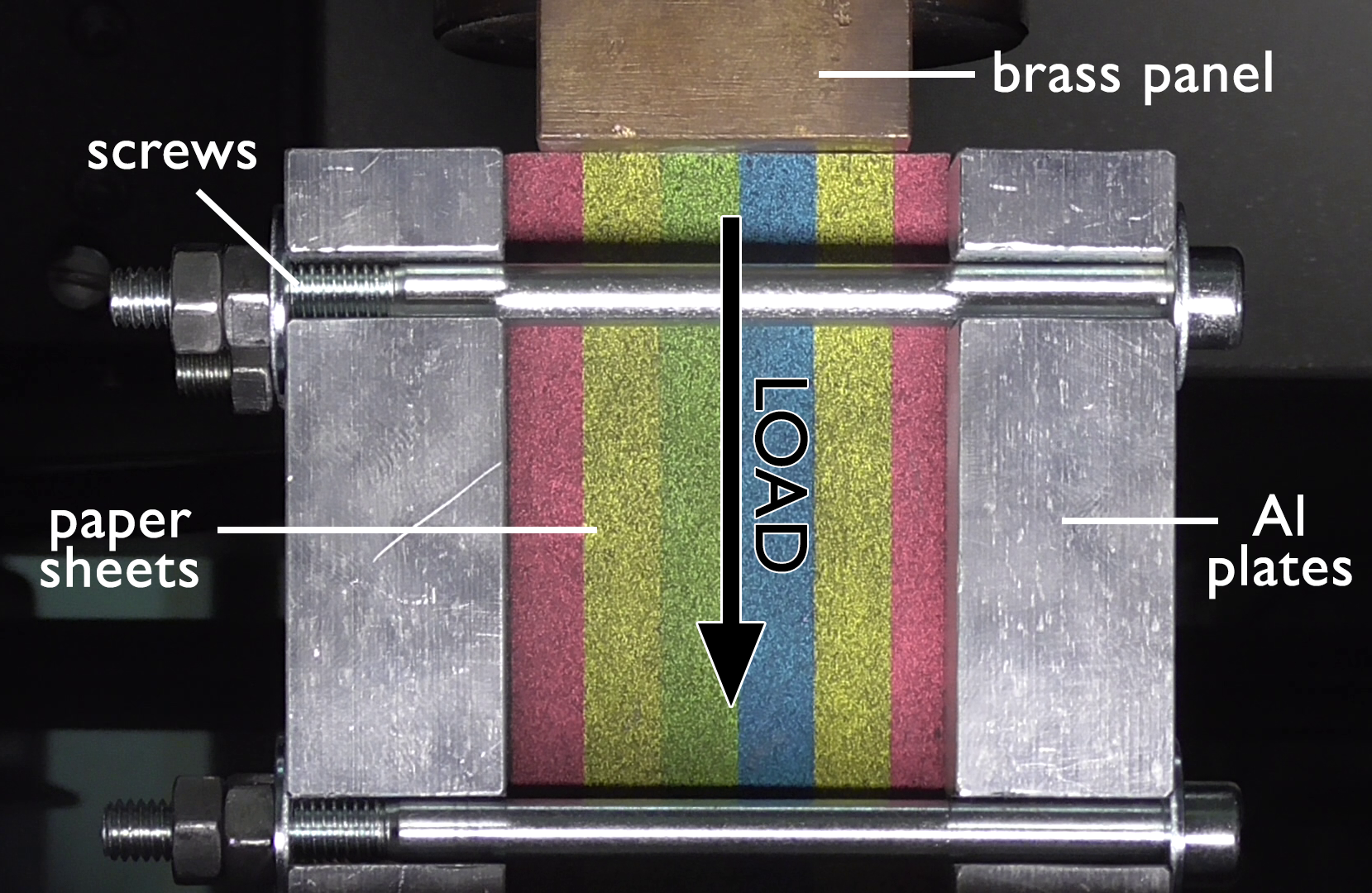}
\caption{Experimental setup of the compression test.}
\label{fig_comprSetup}
\end{center}
\end{figure}

Following the idea of \cite{wadee} we perform a compression experiment on a stack of papers.
Compression tests were performed using a custom-made experimental setup, where blocks of paper sheets were used as testing material; see Fig. \ref{fig_comprSetup}. The dimensions of the block were 75 $\times$ 75 $\times\sim$42 mm$^3$ (height $\times$ width $\times$ thickness) and placed between two thick aluminum plates so that the sheets were parallel with the plates and placed on the steel platen. The aluminum plates were bolted together using four steel screws (near each corner). Thus, the block was confined from the sides and from the bottom, and loading was applied from the top by means of a brass panel slightly longer and narrower than the top face of the block. Three compression tests were carried out using the Instron 5582 universal testing machine with a constant cross-head speed of 0.075 mm$\cdot$s$^{-1}$ giving an initial strain rate of 10$^{-3}$$\cdot$s$^{-1}$.

The video recording was performed using the Panasonic HC-V180EP-K camera with full HD resolution (1920 $\times$ 1080 px$^2$). The front side of the paper block (i.e., perpendicular to the normal of the sheets) facing the camera was sprayed with black acrylic paint in order to create a speckle pattern for the digital image correlation (DIC) analysis. \emph{Ncorr} scripts implemented in MATLAB were used for the computation of DIC with the following parameters: subset size of 15 px and spacing of 5 px. A detailed description of the \emph{ncorr} DIC algorithms is available in Refs. \cite{Blaber2015, ncorr}. 

The obtained compression curves are shown in Fig. \ref{fig_simres3} (left) and exhibit an intermittent plasticity character manifested by means of several pronounced stress drops during loading. These results, supplemented by the DIC data, show that serration corresponds to the gradual formation and propagation of localized ``zig-zag'' deformation bands referred to as kinking, as can be seen in Fig. \ref{fig_simres1} and \ref{fig_simres2}. It can be observed in Supplementary video \cite{figshare} that all load drops correlate with the formation of kink bands. Kinking is a deformation mechanism often observed in layered systems of different scales when compressively loaded parallel to their basal planes \cite{hagihara,plummer}. During loading, the layers first undergo an elastic buckling which is then followed by a nucleation of so-called ripplocation boundaries driven by strain redistribution from high-energy in-plane bonds to low-energy out-of-plane bonds \cite{plummer}.

\section{Numerical verification}

\begin{figure}
\begin{center}
\includegraphics[width = 0.8\textwidth]{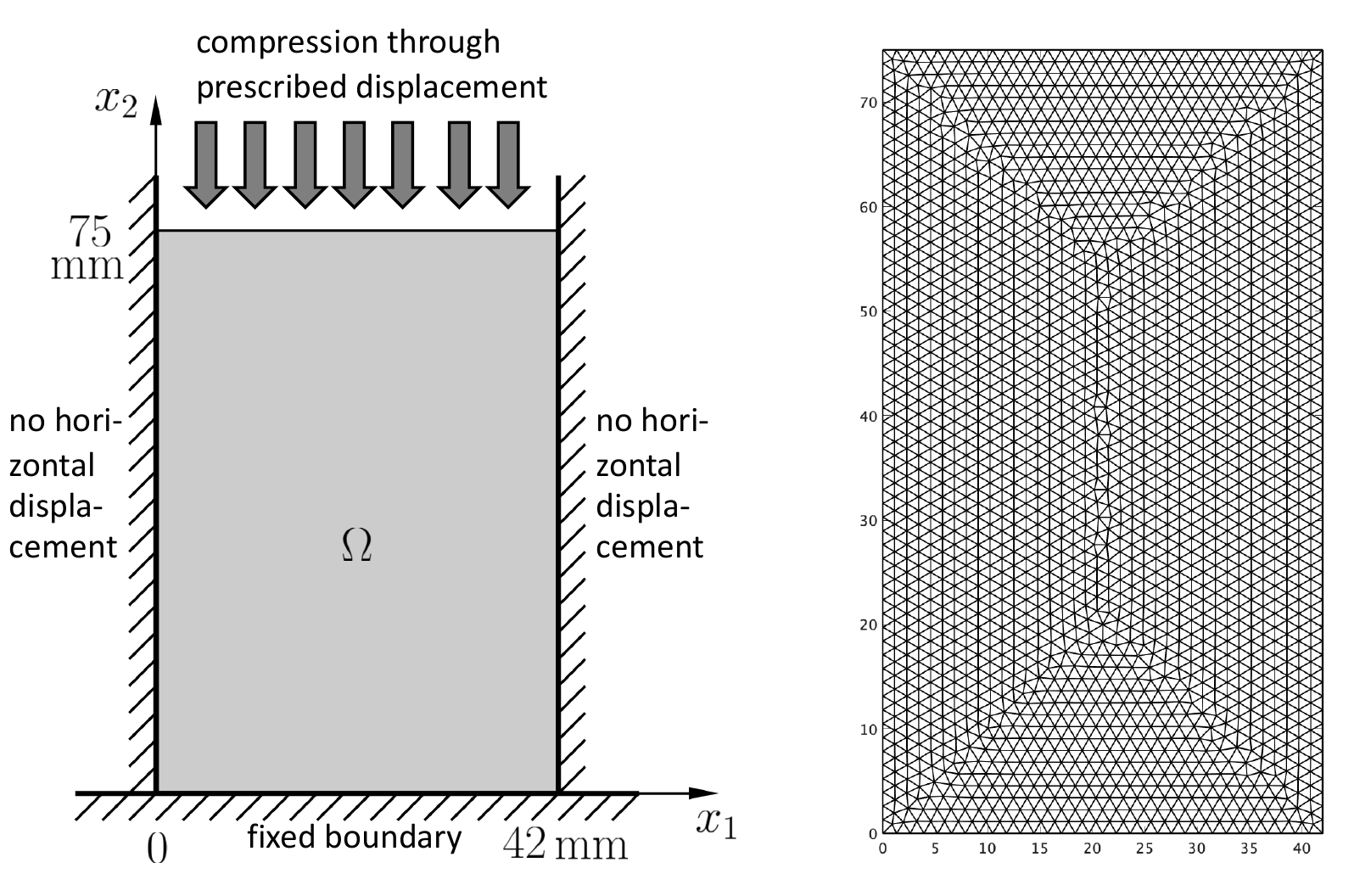}
\caption{Setup of the numerical test and example of a finite element mesh.}
\label{fig_setup}
\end{center}
\end{figure}

We implement a numerical method for the mathematical model to verify that the model is capable of reproducing the main features of the results obtained in the physical experiment. The model was suggested in \cite{conti-dolzmann-kreisbeck}
Since we wish to relate the outputs also to the fully rigid analysis of \cite{conti} and since the precise elastic constants of the material used in experiments are unknown, we abandon the quest for quantitative comparison and adopt the 2D setting in our computations, i.e., $d=2$.
Imitating the setup of the experiment, we take a rectangular region $\Omega$ of material and subject it to loading through a prescribed displacement of the upper edge in the vertical direction.
The lower edge is fixed and the lateral edges are allowed to move only in the vertical direction (see Figure \ref{fig_setup}).
The elasticity of the material is assumed to follow the neo-Hookean model with a single slip system having glide direction $s = (0,1)^T$ and slip-plane normal $m = (1,0)^T$. 
External forces are not considered, i.e., $L=0$,  except for the prescribed time-dependent Dirichlet boundary condition at the upper edge of the domain, from which the energy functional \eqref{energy_gen} takes the specific form
\begin{multline}
\label{num_func}
\mathcal{I} (t,y,\gamma) = \int_{\Omega} \Big\{ C \left( |F^e|^p-d^{p/2} - 2\log (\det F^e) \right) + D (\det F^e -1)^2 + \alpha \operatorname{tr} \left( (F^e)^TF^e M \right) \\
+ \beta | F^p |^r +  \epsilon | \nabla F^p|^{\alpha} \Big\} \, {\rm d}x 
\end{multline}
if $\det F^e>0$ and $+\infty$ otherwise.
We consider $p>2$, $\alpha=r=2$ that implies that $\gamma\in L^{\tilde r}(\O)$ for every $1\le\tilde r<+\infty$, so then we can assume that, in fact, $2\le r<+\infty$ is arbitrary.
Here, $q=(y,\gamma)$ is the state described by deformation $y\in W^{1,p_y}(\Omega;\R^2)$ and slip $\gamma\in W^{1, 2}(\Omega)$, 
coefficients $C,D,\alpha,\beta,\epsilon$ are positive material constants and the elastic and plastic parts of deformation gradient read
\begin{eqnarray*}
&&F^e(y,\gamma) = \nabla y(x) ( I - \gamma(x) s \otimes m ), \qquad F^p(\gamma) = I +\gamma(x) s \otimes m.
\end{eqnarray*}
Furthermore, the term containing $M$ in  \eqref{num_func} is taken from \cite{schroeder} and describes  the transverse isotropy inherent to the stack of paper sheets; here the matrix $M = m\otimes m$ characterizes the material symmetry along the direction $m=(1,0)^T$. 
The integrand of \eqref{num_func} is polyconvex and thus microstructure formation 
will not occur. 
In conjunction with the coercivity of the functional, the existence of a minimizer in $\Y$ for a fixed $\gamma \in \Z$ is guaranteed.

In order to calculate the rate-independent evolution of the system, we employ the dissipation distance between two states with slip magnitudes $\gamma_1$ and $\gamma_2$ defined by
$$ \D(\gamma_1,\gamma_2)=\sigma \int_\Omega |\gamma_1(x)-\gamma_2(x)|\,{\rm d}x\ . $$
In the numerical simulations, we discretize the given time interval $[0,T]$ into $K$ subintervals of equal length $\tau = T/K$  separated by the time points $t^k = k\tau$, $k=0,1, \dots, K$. 
The time-discretization of the rate-independent evolution model then leads to the minimization problem
\begin{equation}
\label{num_minprob}
    q^{k+1} \in \operatorname{Arg} \min_{\widetilde{q}=(\widetilde{y},\widetilde{\gamma})} \left( {\mathcal{D}}(\gamma^k,\widetilde{\gamma}) + \mathcal{I}(t^{k+1},\widetilde{y},\widetilde{\gamma}) \right),
\end{equation}
where the superscript $k+1$ denotes the time instances $ t^{k+1}$ at which the quantity is evaluated, while the dependence of the energy $\mathcal{I}$ on $t^{k+1}$  emphasizes the presence of a time-dependent Dirichlet boundary condition. 
As mentioned above, it can be shown that the time interpolants of its solutions converge to the corresponding energetic solution.

For spatial discretization, the rectangular region is partitioned into a regular triangular mesh and the functionals $\mathcal{I}, \mathcal{D}$ are approximated using the finite element method (FEM) \cite{ciarletfem}. 
In particular, we use linear $\mathcal{P}^1$-elements for both the position vector variable $y(x)$ and the scalar-valued slip variable $\gamma(x)$, obtaining their piece-wise affine approximations $\widehat{y}(x), \widehat{\gamma}(x)$.
Hence, the trial functions for the discretized state $q(x)$ have the form
$$ \widehat{q}(x) = (\widehat{y}_1(x), \widehat{y}_2(x), \widehat{\gamma} (x)) \in V_h \times V_h \times V_h, $$ 
where $V_h$ stands for the subspace of piece-wise affine functions on the triangular partition.
We denote the standard FEM basis functions of $V_h$ by $\{\varphi_i\}_{i=1}^{N}$, where $N$ is the number of nodes in the mesh, i.e., $\varphi_i$ is a hat-like function with value $1$ at the $i$-node and with value $0$ at all other nodes. 
Then we can write
\begin{equation}
\label{num_basisexp}
   \widehat{y}_j(x) = \sum_{i=1}^N a^i_j \varphi_i(x), \; j=1,2, \qquad \widehat{\gamma}(x) = \sum_{i=1}^N b_i \varphi_i(x) 
\end{equation}
for some set of coefficients $\{a^i_1,a^i_2,b_i\}_{i=1}^N \subset \R$.
All coefficients $\{ b_i\}_{i=1}^N$ are free, while some of the values of $a^i_j$ follow the imposed boundary conditions. 
Namely, for nodes $i$ on the lower edge of the domain $a_1^i, a_2^i$ are kept constant (equal to initial values) throughout the computation to realize a fixed boundary, for nodes on the upper edge $a_1^i, a_2^i$ are chosen so that the time-evolving Dirichlet boundary condition is satisfied, and finally, for nodes on the lateral edges, we keep $a_1^i$ equal to their initial values to express the restriction on horizontal displacements.
Hence, denoting the set of indices of free coefficients $a_1^i, a_2^i$ by $I_1, I_2$, respectively, and plugging the expansions \eqref{num_basisexp} into the functional in \eqref{num_minprob}, we obtain a nonlinear functional
\begin{multline}
\label{num_discrmin}
    H^k
    \left(\{a_1^i\}_{i\in I_1}, \{ a_2^i\}_{i\in I_2}, \{b_i\}_{i=1}^N \right) \\
    = {\mathcal{D}}^{\delta} \left( \sum_{i=1}^N b_i^k \varphi_i, \sum_{i=1}^N b_i \varphi_i \right) + \mathcal{I}^{k+1} \left( \sum_{i=1}^N a^i_1 \varphi_i, \sum_{i=1}^N a^i_2 \varphi_i, \sum_{i=1}^N b_i \varphi_i  \right)
\end{multline}
to be minimized 
over the free variables $\{a_1^i\}_{i\in I_1}, \{ a_2^i\}_{i\in I_2}, \{b_i\}_{i=1}^N$ at each time level.
We remark that the functional $\mathcal{D}$ has been replaced in the numerical implementation by its smoothed version
\begin{equation}
    \label{num_regul}
\mathcal{D}^{\delta} (\gamma_1,\gamma_2)=\sigma \int_\Omega \sqrt{\delta^2 + |\gamma_1(x)-\gamma_2(x)|^2} \,{\rm d}x\ ,
\end{equation}
where $\delta$ is a small positive parameter. 

After initialization of the coefficient $a_j^i, b_i$ to represent the undeformed configuration, i.e., setting $(a_1^i,a_2^i)$ equal to the $(x,y)$-coordinates of $i$-th mesh node and $b^i=0$ for all nodes, the numerical computation at further time steps proceeds as follows. 
At the beginning of each time step, the coefficients $a_j^i$ in \eqref{num_basisexp} corresponding to the mesh nodes located on the upper edge of the domain are updated to represent the prescribed displacement at the current time $t^k$.
Subsequently, the resulting function $H^k$ from \eqref{num_discrmin} is minimized over the free variables $\{a_1^i\}_{i\in I_1}, \{ a_2^i\}_{i\in I_2}, \{b_i\}_{i=1}^N$, and the result is fed as input to the next time step.
The minimization is carried out using the function {\sl fminunc} from MATLAB's Optimization toolbox \cite{matlab}.
This function is based on a trust-region optimization algorithm \cite{conn}, which requires as input the function $H^k$, its gradient, initial guess for solution, and optionally also the sparsity pattern of the Hessian of $H^k$. 
The function $H^k$ and the finite-difference approximation of its gradient are evaluated using the vectorized algorithm of \cite{mosval}, which produces a significant speed-up of the computation. 
The initial guess is taken as the solution in the previous time step for the elastic part $a_j^i$, and as zero for the plastic part $b_i$.
In the minimization procedure, it might be physically meaningful to distinguish between two approaches: to minimize over both elastic and plastic variables at once, or to alternately minimize over the elastic or plastic part while keeping the other part temporarily fixed.
However, in this case, no difference was observed in the impact of the order of minimization on the numerical results.

In the specific simulation, we consider the deformation in the time interval $t\in [0,100]$ (in seconds) of the rectangular region $\Omega = (0,42)\times(0,75)$ (in millimeters), whose upper edge is compressed with constant speed according to $y(t,x_1,x_2) = (x_1,75-0.18t)$.
Time is discretized into $K=76$ steps of length $\tau=1.3125$, and the undeformed region $\Omega$ is divided into a triangular mesh with 4184 elements and 2182 nodes (see Figure \ref{fig_setup}).
The values of the remaining parameters are selected as follows:
material constants in the energy functional \eqref{num_func} are $C=0.6$ GPa, $D =0.2$ GPa, $\alpha = 0.1$ GPa, $\beta = 20$ kPa, $\epsilon =500$ N, $\sigma = 1$ kPa, 
the regularization parameter in \eqref{num_regul} is $\delta = 10^{-5}$, 
the perturbation to calculate a finite difference approximation of the gradient of $H^k$ is $10^{-8}$, 
the infinite value for $\det F^e \leq 0$ in \eqref{num_func} is realized by penalty $10^6$ on individual triangles,
integration over triangles uses the second-order quadrature rule,
and finally stopping criteria in the trust-region minimization algorithm are TolX $=10^{-10}$ for step and TolFun $=10^{-4}$ for function values.
The elastic constants are chosen so that their order of magnitude corresponds to the experimentally measured values reported in \cite{schulgasser, yokoyama}.

\begin{figure}[h!]
\centering
\includegraphics[width = 0.8\textwidth]{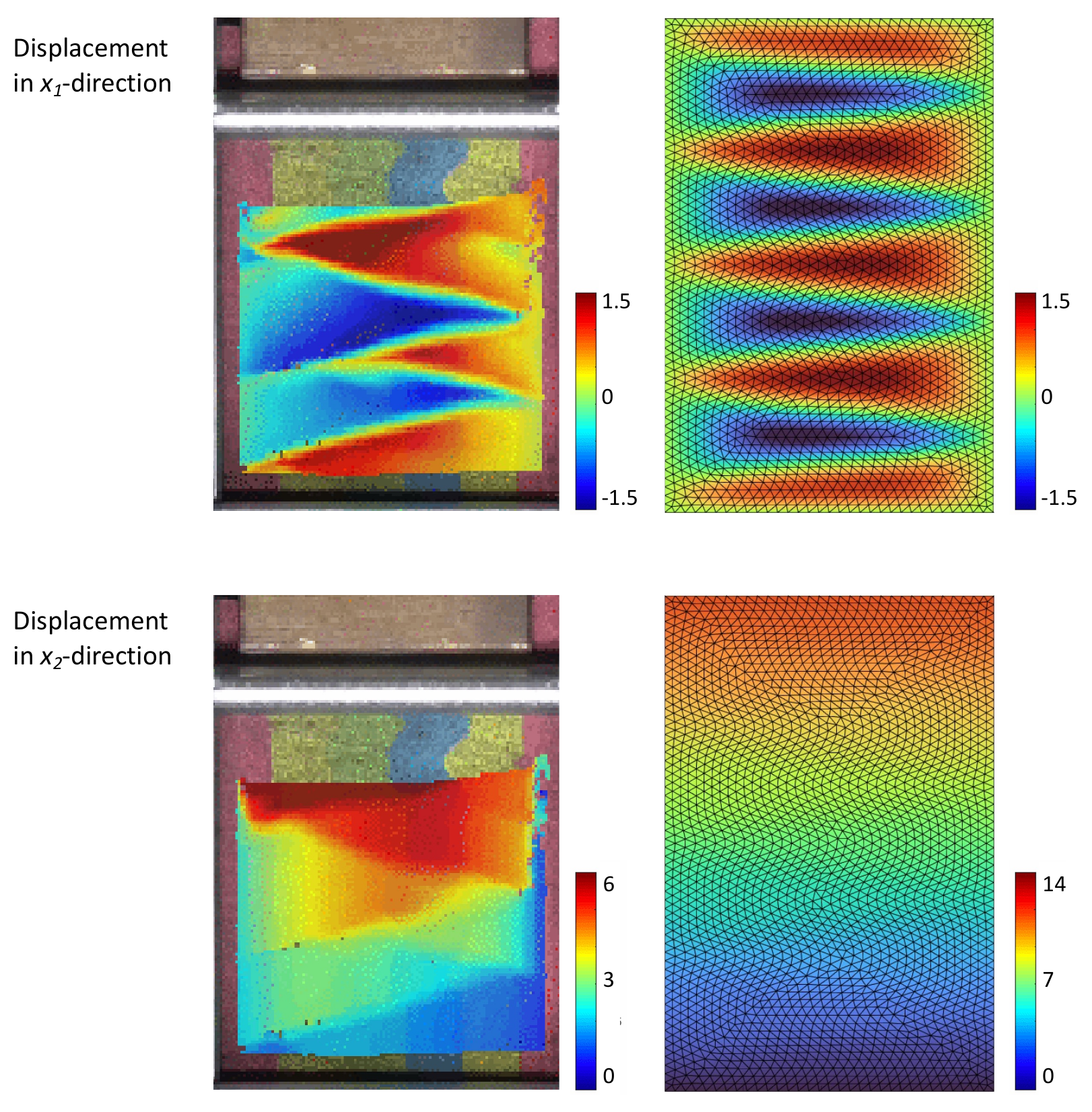}
\caption{Comparison of displacements in $x_1$- and $x_2$-directions obtained in the experiment and numerical simulation. The images in the left column show experimental observations at the time when the vertical displacement of upper boundary is around 12 mm and when the kink-bands reach the lower boundary. The right images show numerical results at the same vertical displacement of upper boundary by 12 mm, a little after the kink-band structure appears in the simulation.}
\label{fig_simres1}
\end{figure}

\begin{figure}[h!]
\begin{center}
\includegraphics[width = 0.8\textwidth]{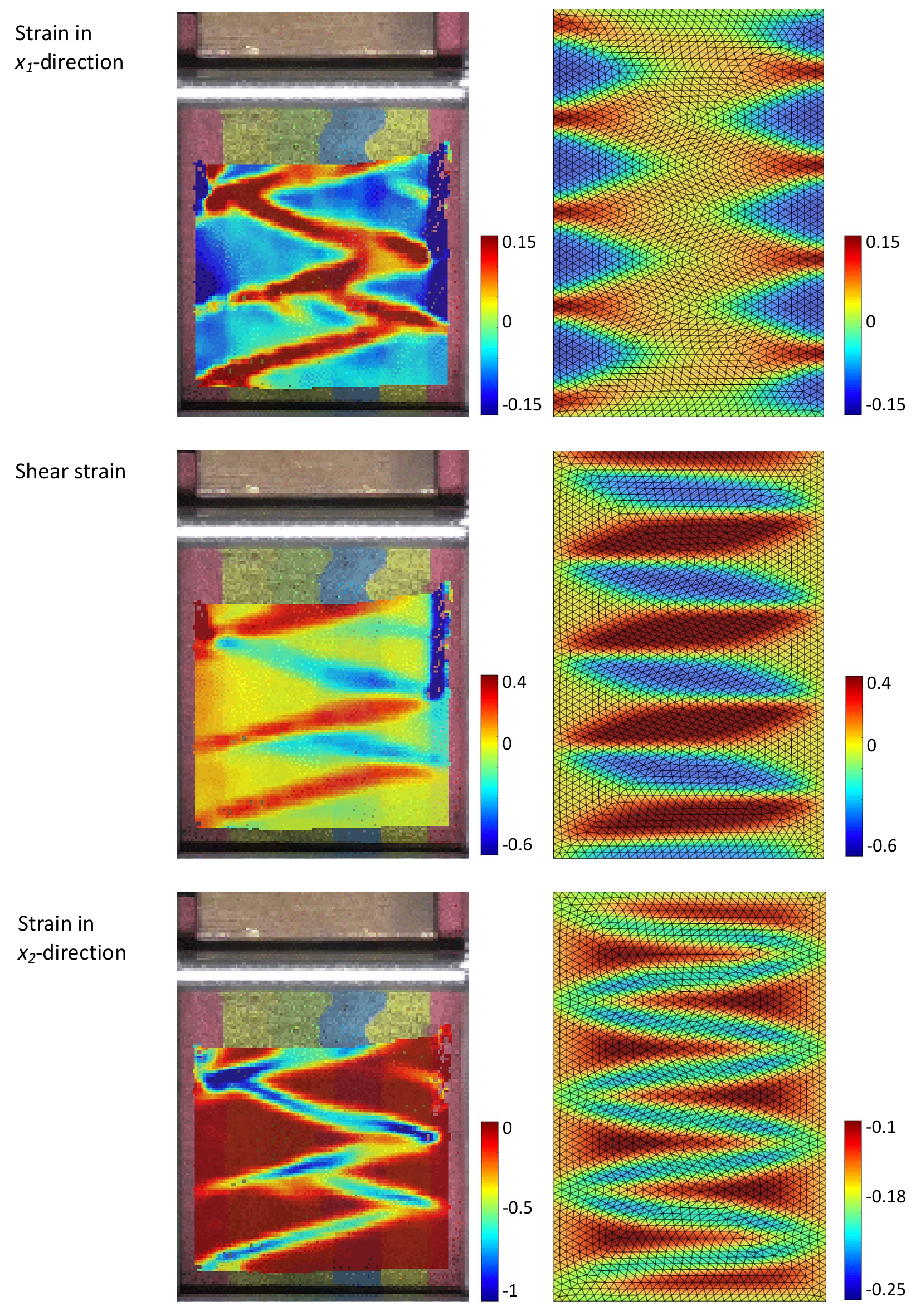}
\caption{Comparison of strains obtained in the experiment (left column) and in numerical simulation (right column). The images were taken at the same load conditions as in Figure \ref{fig_simres1}.}
\label{fig_simres2}
\end{center}
\end{figure}

\begin{figure}[h!]
\begin{center}
\includegraphics[width = 0.95\textwidth]{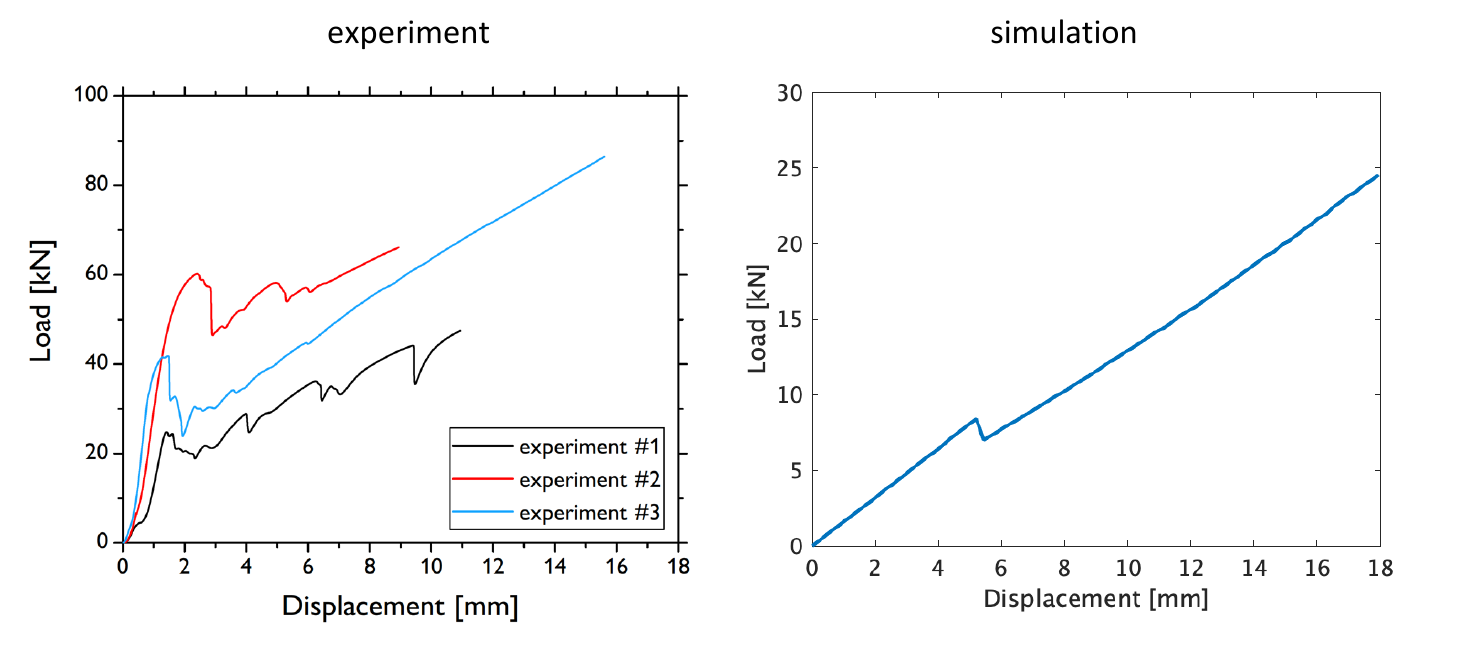}
\caption{Comparison of stress-strain curves obtained in the experiment and in numerical simulation. Curves for three independent experiments are shown, where the lowest curve (experiment \#1) corresponds to images displayed in Figures \ref{fig_simres1} and \ref{fig_simres2}. In experiment, kinks form gradually one by one due to finite speed of sound, while in simulation they form all at once, which leads to a later appearance of kinks and hence a slightly later drop in the load.}
\label{fig_simres3}
\end{center}
\end{figure}

The results obtained by implementing the above described model are summarized in Figures \ref{fig_simres1} (for displacements), \ref{fig_simres2} (for strains) and \ref{fig_simres3} (for displacement-load curves), together with the results of the experimental analysis.
Despite a simple mathematical model, surprisingly good qualitative and even quantitative agreement is reached.
In particular, the simulation succeeds in reproducing the formation of wedge-like kink-bands of the same type as in experiment, and the displacements and strains show identical deformation patterns.
The lack of an accurate quantitative match is mainly due to the fact that only a two-dimensional model is used and the detailed elastic properties of the paper sheets used in the experiment are unknown.

We remark that the color-bar range for the displacement in $x_2$-direction in Figure \ref{fig_simres1} is different for experiment and simulation mainly due to the fact that only a part of the material region is analyzed owing to limitations of the experimental setup.
Even if the region is adjusted, a complete match is not gained, as can be predicted from the difference in $x_2$-strain in Figure \ref{fig_simres2} -- notice the different range of color-bars for experiment and simulation.
This discrepancy is thought to be caused mainly by the fact that in experiment the stack of papers is not compressed in a thin region at the left and right boundaries (this is visible in the upper left and upper right corners of the figure), which allows for a relaxation of stresses near the lateral boundaries.

Regarding the dynamics, in the real phenomenon kink-bands are formed one by one in sequence, as can be confirmed in the repository \cite{figshare}.
On the other hand, the mathematical model is quasi-stationary, and thus all the kinks appear at once as soon as sufficient energy is stored to allow for the plastic deformation.
Nevertheless, these issues can be addressed by a straightforward refinement of the basic model, and therefore, the results show that the proposed mathematical model is capable of capturing the important features of the deformation mechanism. 
Moreover, the results indicate that the rate-independent evolution provides an "elastically regularized" approximation for the completely rigid problem of \cite{conti}.


Movies of the experiments and numerical simulation are available at 
\href{https://figshare.com/projects/Elastoplastic_deformations_of_layered_structure/132212}{Figshare} \cite{figshare}.

\section*{Acknowledgements}
Martin Kru\v{z}\'{i}k and Jan Valdman were supported by the GA\v{C}R project 21-06569K.
Martin Kru\v{z}\'{i}k also thanks the ESI Vienna for its hospitality during his stay in January-February 2022.
The research of Karel \v{S}vadlenka was supported by JSPS Kakenhi Grant numbers 19K03634 and 18H05481.
Support in the framework of Visegrad Group (V4)-Japan Joint Research Program -- Advanced Materials under grant No. 8F21011 is gratefully acknowledged by Daria Drozdenko and Kristi\'{a}n M\'{a}this.

\end{document}